\documentclass[10pt, reqno]{amsart}

          \usepackage{fullpage}
          \usepackage{amssymb}
          \usepackage{amsmath}
          \usepackage{amsthm}
          \usepackage{amsfonts}
          \usepackage{enumerate}
          \usepackage{url}
          \usepackage{color}
          \usepackage{graphicx}
          \usepackage[dvipsnames]{xcolor}

\usepackage[numbers,sort&compress]{natbib}

\newtheorem{theorem}{Theorem}

\newtheorem{lemma}{Lemma}

\theoremstyle{definition}

\newcommand{\comment}[1]{}

\newcommand{\ind}{{\bf 1}}

\def\inddd#1{{\ind}_{\left\{#1\right\}}}

\newcommand{\proba}{\mathbb P}
\newcommand{\PR}{\mathbb P}

\newcommand{\esp}{{\mathbb E}}

\newcommand{\inv}{^{-1}}
\newcommand{\cov}{{\rm{Cov}}}
\newcommand{\var}{{\rm{Var}}}

\newcommand{\eqnh}{\begin{eqnarray*}}
\newcommand{\eqne}{\end{eqnarray*}}
\newcommand{\eqnhn}{\begin{eqnarray}}
\newcommand{\eqnen}{\end{eqnarray}}
\newcommand{\equh}{\begin{equation}}
\newcommand{\eque}{\end{equation}}

\def\summ#1#2#3{\sum_{#1 = #2}^{#3}}

\def\sif#1#2{\sum_{#1=#2}^\infty}

\newcommand{\eqd}{\stackrel{d}{=}}

\def\nn#1{{\left\|#1\right\|}}

\def\ccbb#1{\left\{#1\right\}}

\def\pp#1{\left(#1\right)}

\def\bb#1{\left[#1\right]}

\def\mmid{\;\middle\vert\;}

\def\floor#1{\left\lfloor #1 \right\rfloor}

\def\ceil#1{\left\lceil #1 \right\rceil}

\def\vv#1{{\boldsymbol #1}}

\def\qmand{\quad\mbox{ and }\quad}

\def\mfa{\mbox{ for all }}

\def\mmas{\mbox{ as }}

\def\wt#1{\widetilde{#1}}
\def\wb#1{\overline{#1}}
\def\what#1{\widehat{#1}}

\def\limn{\lim_{n\to\infty}}

\def\weakto{\Rightarrow}

\def\R{{\mathbb R}}

\def\N{{\mathbb N}}

\def\B{{\mathbb B}}

\def\G{{\mathbb G}}

\def\calK{\mathcal K}
\def\calL{\mathcal L}

\def\calN{\mathcal N}

\usepackage{epsfig}
\usepackage{ifthen}
\theoremstyle{definition}

\numberwithin{equation}{section}

\def\Q{{\mathbb Q}}
\def\F{{\mathbb F}}
\title{Remarks on power-law random graphs}

\author{Mei Yin}
\address{Department of Mathematics, University of Denver, Denver, CO 80208}
\email{mei.yin@du.edu}

\begin{document}

\maketitle

The theory of graphons is an important tool in understanding properties of large networks. We investigate a power-law random graph model and cast it in the graphon framework. The distinctively different structures of the limit graph are explored in detail in the sub-critical and super-critical regimes. In the sub-critical regime, the graph is empty with high probability, and in the rare event that it is non-empty, it consists of a single edge. Contrarily, in the super-critical regime, a non-trivial random graph exists in the limit, and it serves as an uncovered boundary case between different types of graph convergence.

\vspace{.2cm}


\vspace{.2cm}


\section{Introduction}\label{intro}
The advent of large networks in the last decades has posed exciting challenges for researchers. Many interesting questions have been asked, ranging from a notion of limit distribution for sequences of graphs, to the understanding of local and global characteristics of large graphs, as well as the existence of effective algorithms to generate graphs with desired properties, and many more. Following the earlier work of Aldous \cite{Aldous} and Hoover \cite{Hoover}, Lov\'{a}sz and coauthors (V.T. S\'{o}s, B. Szegedy, C. Borgs, J. Chayes, K. Vesztergombi, $\dots$) have constructed an elegant theory of graph limits in a sequence of papers \cite{borgs08convergent} \cite{borgs12convergent} \cite{LS} to address these questions. The graph limit theory sheds light on various topics such as graph testing and extremal graph theory, and has found applications in statistics and related areas (see for instance Bhattacharya et al. \cite{B}, Bickel and Chen \cite{BC}, Chatterjee et al. \cite{CDS}, and Mukherjee and Xu \cite{MX}). For a comprehensive study on the theory of graph limits, we refer to van der Hofstad \cite{Hofstad} and Lov\'{a}sz \cite{Lov}.

We present some basics of the graph limit theory in the context of dense graphs (number of edges comparable to the square of number of vertices), which was the initial setting for this theory. Let $\mathcal{W}$ be the space of all symmetric, measurable functions $W : [0,1]^{2} \to [0,1]$ (referred to as the \emph{graph limit space} or \emph{graphon space}). Any simple graph $G_n$, irrespective of the number of vertices $n$, may be represented as an element $W_n \in \mathcal{W}$, such that
	\begin{equation}
		W_n(x,y) = \left\{ \begin{array}{ll} 1 & \text{if }\left\{\left\lceil nx \right\rceil, \left\lceil ny \right\rceil\right\} \text{is an edge in } G_n, \\ 0 & \text{otherwise.} \\  \end{array} \right.
	\end{equation}
A sequence of graphs $\{G_n\}_{n\geq 1}$ is said to converge to a function $W \in \mathcal{W}$ if for every finite simple graph $H$ with vertex set $V(H)=[k]=\{1,...,k\}$ and edge set $E(H)$,
\begin{equation}\label{sub}
\lim_{n \to \infty} t(H, W_n)=t(H, W),
\end{equation}
where by construction,
\begin{equation}
t(H, W_n)=\frac{|\text{hom}(H,
G_n)|}{|V(G_n)|^{|V(H)|}}
\end{equation}
equals the density of graph homomorphisms from $H$ to $G_n$, and
\begin{equation}
\label{tt} t(H, W)=\int_{[0,1]^k}\prod_{\{i,j\}\in E(H)}W(x_i,
x_j)dx_1\cdots dx_k.
\end{equation}
There are generalized versions of this type of convergence for weighted graphs, but the principal ideas are the same. Every function in $\mathcal{W}$ is the limit of a certain convergent graph sequence \cite{LS}. Intuitively, the interval $[0,1]$ represents a continuum of vertices, and $W(x,y)$ denotes the probability of putting an edge between $x$ and $y$. For example, for the Erd\"{o}s-R\'{e}nyi random graph $G(n,\rho)$, the associated limiting graphon is represented by the function that is identically equal to $\rho$ on $[0,1]^2$.

The graphon interpretation enables us to capture the notion of convergence in terms of subgraph densities (\ref{sub}) by an explicit metric on $\mathcal{W}$, the \emph{cut distance}:
\begin{equation}
d_{\square}(U, V)=\sup_{S, T \subseteq [0,1]}\left|\int_{S\times
T}\left(U(x, y)-V(x, y)\right)dx\,dy\right|
\end{equation}
for $U, V \in \mathcal{W}$. A non-trivial complication is that the topology induced by the cut metric is well defined only up to
measure preserving transformations of $[0,1]$ (and up to sets of Lebesgue measure zero), which in the context of finite graphs may
be thought of as vertex relabeling. To tackle this issue, an equivalence relation $\sim$ is introduced in $\mathcal{W}$. We say
that $U\sim V$ if $U(x, y)=V_{\sigma}(x, y):=V(\sigma x, \sigma y)$ for some measure preserving bijection $\sigma$ of $[0,1]$. Let
$\tilde{U}$ (referred to as a \emph{reduced graphon} or \emph{unlabeled graphon}) denote the closure of the orbit $\{U_\sigma\}$ in $(\mathcal{W}, d_{\square})$. Since $d_{\square}$ is invariant under $\sigma$, one can then define on the resulting quotient space $\tilde{\mathcal{W}}$ the natural distance $\delta_{\square}$ by $\delta_{\square}(\tilde{U}, \tilde{V})=\inf_{\sigma_1, \sigma_2}d_{\square}(U_{\sigma_1},
V_{\sigma_2})$, where the infimum ranges over all measure preserving bijections $\sigma_1$ and $\sigma_2$, making
$(\tilde{\mathcal{W}}, \delta_{\square})$ a metric space. With some abuse of notation we also refer to $\delta_{\square}$ as
the \emph{cut distance}. The space $(\tilde{\mathcal{W}}, \delta_{\square})$ enjoys many nice properties. For example, it is a compact space and homomorphism densities $t(H, \cdot)$ are continuous functions on it.

In addition to developing a standard theory of limits for sequences of dense graphs, serious efforts have been made at formulating parallel results for sparse graphs. Unlike dense graphs, sparse graphs display vastly different edge densities. Nevertheless, all simple sparse graphs converge to the zero graphon in the classical theory of graphons. To take into account the wide-ranging edge densities of sparse graphs, a renormalization procedure on their associated graphons has been introduced. The renormalization is executed in two ways, either by rescaling the height of the graphon or by stretching the domain on which it is defined.

The first renormalization approach was introduced in Bollob\'{a}s and Riordan \cite{BR} and Borgs et al. \cite{borgs19Lp}, and is characterized by the {\em rescaled cut metric} $\delta^r_\square$. The graphon representation $W_n$ of $G_n$ is rescaled to $W_n^r$ by
\equh\label{eq:rescaled}
W_n^r(x, y)=\nn {W_n}_1\inv W_n(x, y),
\eque
i.e., divide the weight of each edge by the edge density $\nn {W_n}_1$. The rescaled limit graphon $W^r$ is in $L^1([0,1]^2)$, and takes values in $[0,\infty)$. Note that for sparse graphs this scaling of the edges is appropriate, as it produces a constant order for $\nn {W_n^r}_1$. If we rescale $W_n$ so that $\nn {W_n^r}_1=o(1)$ instead, then automatically its cut norm (which is upper bounded by the $L^1$-norm) is negligible. As for dense graphs, there are generalized versions of this type of convergence for weighted sparse graphs, as well as to $L^p$ graphons.

The second renormalization approach was introduced in Borgs et al. \cite{borgs17sparse}, and is characterized by the {\em stretched cut metric} $\delta^s_\square$. This time, the graphon representation $W_n$ of $G_n$ is stretched to $W_n^s$ by
\equh\label{eq:stretched}
W_n^s(x,y)= W_n\pp{\nn {W_n}_1^{1/2}x, \nn {W_n}_1^{1/2}y},
\eque
i.e., multiply the input arguments by the square root of the edge density $\nn {W_n}_1$ for valid $x, y$ values. Equivalently, the stretching may be interpreted as rescaling the measure of the underlying measure space. The stretched limit graphon $W^s$ is in $L^1([0,\infty)^2)$, and takes values in $[0, 1]$. Under this rescaling of the measure perspective, graphons on $\sigma$-finite measure spaces may be considered as limiting objects for sequences of sparse graphs, similarly as graphons on probability spaces are considered as limits of dense graphs. Again there are further generalizations for this type of convergence.

While introducing the rescaled convergence mode in \cite{borgs19Lp}, Borgs et al. presented a motivating example. Consider a discrete graph $G_n$ of $n$ vertices numbered $1$ through $n$. Connect vertices $i,j$ with probability
\begin{equation}
p_n(i, j)=\min\ccbb{1,n^\beta/(ij)^{1/\alpha}} = \min\ccbb{1,n^{\beta-2/\alpha}(i/n)^{-1/\alpha}(j/n)^{-1/\alpha}},
\end{equation}
where $\alpha>1$ and $\beta \in (0, 2/\alpha)$ are parameters. In other words, the edge connection probability between vertices $i,j$ behaves like $(ij)^{-1/\alpha}$, but boosted by a factor of $n^\beta$ in case it becomes too small. This configuration model is one of the simplest ways to get a power law degree distribution, as the expected degree of vertex $i$ scales according to an inverse power law in $i$ with exponent $1/\alpha$. The parameter range on $\alpha$ and $\beta$ is taken for practical considerations: $\alpha>1$ avoids having almost all the edges of the graph between a sub-linear number of vertices, and $\beta\in (0,2/\alpha)$ ensures that the cut-off from taking the minimum with $1$ affects only a negligible fraction of the edges. Let $\{E_{i,j}\}_{1\le i<j\le n}$ be Bernoulli random variables with parameter $p_n(i,j)$ and set $E_{i,j}=E_{j,i}$. For $x,y\in(0,1)$, let
\begin{equation}
W_n^r(x,y)= \frac1{n^{\beta-2/\alpha}}E_{\ceil{xn},\ceil{yn}}
\end{equation}
be the empirical graphon rescaled by the expected edge density of the graph. The result in \cite{borgs19Lp} (for more details see \cite[Example 3.3.3]{borgs18Lp}) says that
\begin{equation}
\limn W_n^r(x,y) =(1-1/\alpha)^2 (xy)^{-1/\alpha}:=W^r(x, y).
\end{equation}
The convergence is with respect to the cut metric, and the limit graphon $W^r(x,y)$ lies in $L^p([0,1]^2)$ for any $p<\alpha$.

In this paper we will examine a graph model that is closely related to the motivating example discussed above. There are $n$ vertices numbered $1$ through $n$. We adapt the edge connection probability for vertices $i,j$ in two steps as shown below:
\begin{equation}
\min\ccbb{1,n^{\beta-2/\alpha}(i/n)^{-1/\alpha}(j/n)^{-1/\alpha}} \rightarrow \min\ccbb{1,\frac{X_{i}X_{j}}{a_n}} \rightarrow \inddd{\frac{X_{i}X_{j}}{a_n}>1},
\end{equation}
where $a_n = n^{-\beta+2/\alpha}$, $X_i\eqd U_i^{-1/\alpha}$, and $U_i$ are i.i.d.~$(0,1)$-uniform random variables. The first step in the adaptation continualizes the discrete normalized vertex labels into a uniform measure, and implicitly relabels the vertices $1,\dots, n$ using the order statistics of their associated random variables $X_1,\dots,X_n$, the latter not having a real impact on the structure of the graph. The second step in the adaptation is more significant. For every edge that our modified model connects, the original configuration construction in \cite{borgs19Lp} connects them too. Call these ``hard edges''. However, the original construction is not that strict with those edges that we drop. Instead they choose whether to connect them or not depending on a Bernoulli sampling probability between $[0,1]$. Call these ``Bernoulli edges''. It might therefore be apt to refer to the model in \cite{borgs19Lp} as \textit{a power-law random graph with Bernoulli edges} and our adapted model as \textit{a power-law random graph without Bernoulli edges}. Also note that the parameter range investigated in \cite{borgs19Lp} translates to $\alpha>1$ and $a_n\ll n^{2/\alpha}$ in our setting.

At the critical regime ($a_n \sim n^{2/\alpha}$), a realization of this adapted model exhibits a small clique and large numbers of follower vertices asymptotically. The limit structure of the model away from criticality on the other hand is less understood, and will be the central focus of this work. We start with some straightforward calculations in Section \ref{sec:toy} and make some quick observations. Section \ref{sec:sub_critical} studies the sub-critical regime ($a_n\gg n^{2/\alpha}$). We show that although there is no graph in the limit, in the rare event that we do see a non-empty graph, typically it consists of exactly one edge. Section \ref{sec:super_critical} studies the super-critical regime ($a_n\ll n^{2/\alpha}$). We show that unlike the original model in \cite{borgs19Lp}, universality emerges in the limit graphon of the adapted model. After proper scaling, the parameter influence on the relation between number of vertices/edges disappears asymptotically. The qualitative difference between the limit graph structures in the original model vs. the adapted model is essentially due to the presence of Bernoulli edges, as the number of those edges is of larger order than the number of deterministic ones.

Putting our investigation of the asymptotic graph structure into the context of random graph limits, we will see that in the super-critical regime where there is a non-trivial random graph in the limit, our adapted model serves as an uncovered boundary case between different types of graph convergence. In contrast to the power-law random graph with Bernoulli edges analyzed in \cite{borgs19Lp} whose graphon representation converges in the rescaled cut metric, the graphon representation of our power-law random graph without Bernoulli edges converges under a modified stretched convergence mode. See Theorem \ref{thm:1} and the accompanying implications for details.

\section{First estimates}\label{sec:toy}
For $X_i\eqd U_i^{-1/\alpha}$ where $U_i$ are i.i.d.~$(0,1)$-uniform, one could take the i.i.d.~$X_i$ to have pdf $\alpha x^{-\alpha-1}dx\inddd{x\ge 1}$ to make all calculations explicit. Given a realization $X_1, \dots, X_n$ of vertex values and a chosen normalization $a_n$, we group the non-isolated vertices of the graph into two parts depending on whether $X_i>\sqrt{a_n}$ or $X_i \leq \sqrt{a_n}$, respectively referred to as ``clique'' and ``followers''. Since two vertices only get connected when the product of their vertex values exceeds $a_n$, a split graph is produced, as vertices are all connected within the clique and form a complete subgraph, while follower vertices can only be connected to clique vertices but not to themselves. Let us label the vertices according to decreasing $X_i$ vertex values, where $Y_i=X_{i: n}$ is the $i$th largest order statistic of i.i.d.~random variables $X_1, \dots, X_n$. Assume that the clique consists of $K_{n, 0}$ vertices indexed by $1,2,\dots, K_{n,0}$, with vertex values $Y_1\geq \cdots \geq Y_{K_{n,0}}>\sqrt{a_n}$. If $K_{n,0}=0$ then there is no graph and the structure is trivial. So suppose $K_{n,0}>0$. For each $j=1,\dots,K_{n,0}$, add in addition $K_{n,j}$ vertices that each connects to clique vertices indexed by $1,2,\dots,K_{n,0}+1-j$, but not clique vertices indexed by $K_{n,0}+2-j,\dots,K_{n,0}$ or all the other existing vertices. One might for the sake of simplicity set $K_{n,j} = 0$ for all $j>K_{n,0}$. Then, $\vv K_n = \{K_{n,j}\}_{j\in\N_0}$ determines the structure of the random graph completely. For example, the total number of vertices is $|V_n| = \sif j0K_{n,j}$ and the total number of edges is
\begin{equation}
|E_n| = \binom {K_{n,0}}2+ \summ j1{K_{n,0}}(K_{n,0}+1-j)K_{n,j}.
\end{equation}

We further compute the edge connection probability,
\begin{align}
\proba(X_1X_2>a_n)&=\int_1^{a_n} \alpha x_1^{-\alpha-1} dx_1 \int_{\frac{a_n}{x_1}}^\infty \alpha x_2^{-\alpha-1} dx_2+\int_{a_n}^{\infty} \alpha x_1^{-\alpha-1} dx_1 \int_1^\infty \alpha x_2^{-\alpha-1} dx_2
\notag\\&=\alpha a_n^{-\alpha}\log a_n+a_n^{-\alpha} \sim \alpha a_n^{-\alpha}\log a_n,
\end{align}
where ``$\sim$'' indicates that the ratio between the two quantities converges to $1$. At the critical regime $a_n \sim n^{2/\alpha}$, this gives the expected total number of edges as
\begin{equation}
\esp |E_n|=\binom{n}{2}\proba(X_1X_2>a_n) \sim \log n.
\end{equation}
Despite the logarithmic growth of edges, the average size of the clique stays at $1$ asymptotically,
\equh\label{eq:kappa_n,0}
\esp K_{n,0}=n\proba(X_1>\sqrt{a_n}) = \frac{n}{a_n^{\alpha/2}} \sim 1.
\eque
Since $K_{n,0} = \summ i1n \inddd{X_i/\sqrt{a_n}>1}$, it is Binomial distributed with parameter $(n, 1/n)$. The Binomial distribution falls off fast as one moves away from the mean, with $K_{n, 0}=1$ being most probable when a non-trivial graph is produced. Consequently, the average random graph at criticality has a small clique and large numbers of follower vertices.

\begin{figure}
\includegraphics[width = .45\textwidth]
{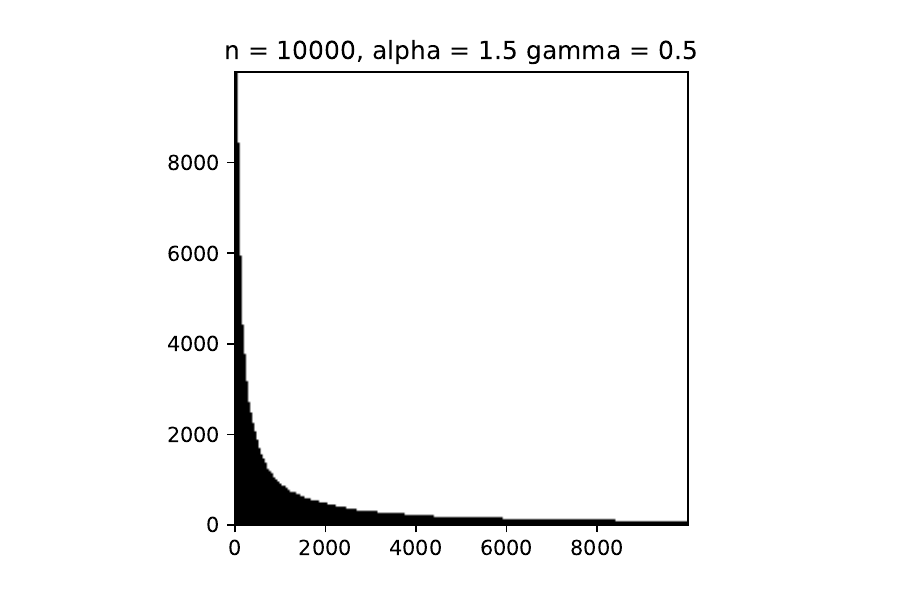}
\hfill
\includegraphics[width = .45\textwidth]
{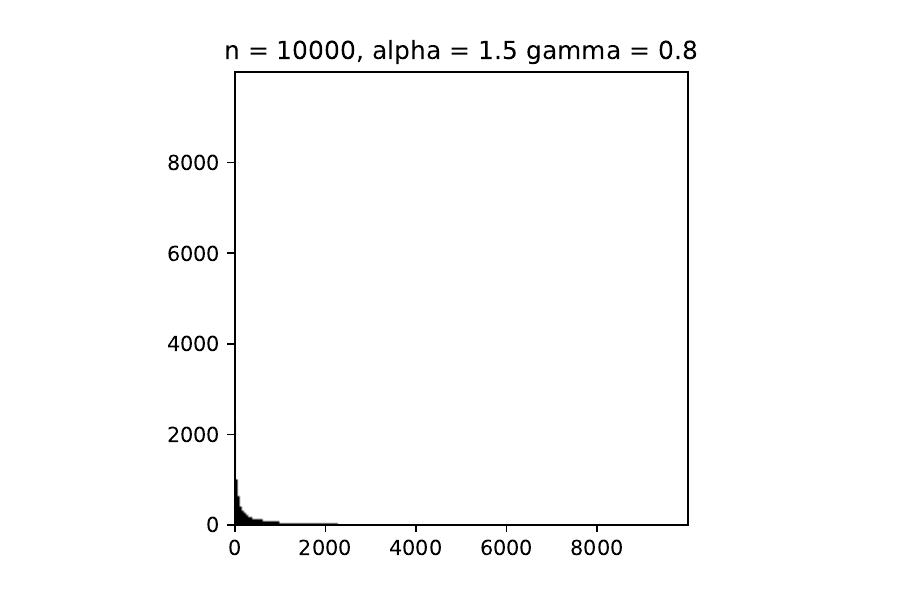}
\includegraphics[width = .45\textwidth]
{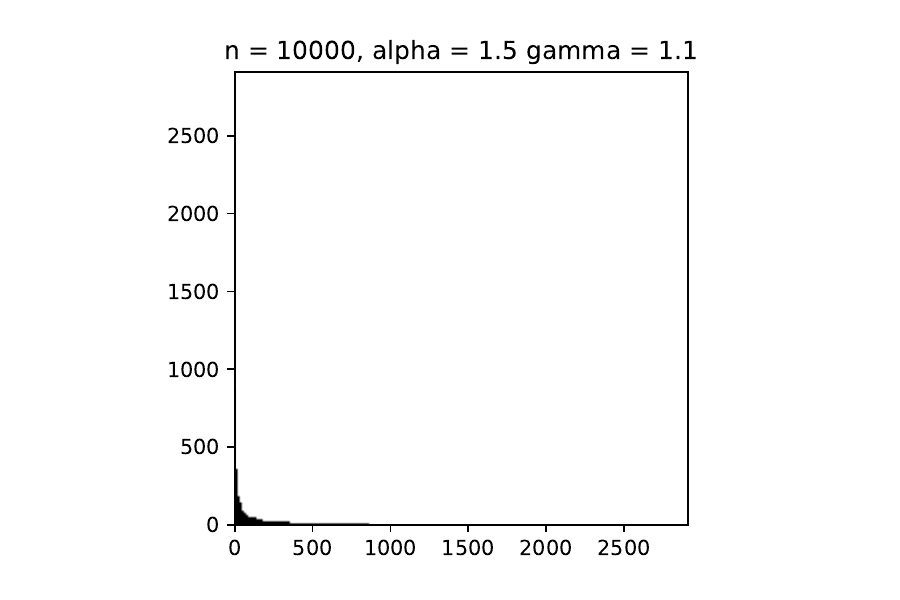}
\hfill
\includegraphics[width = .45\textwidth]
{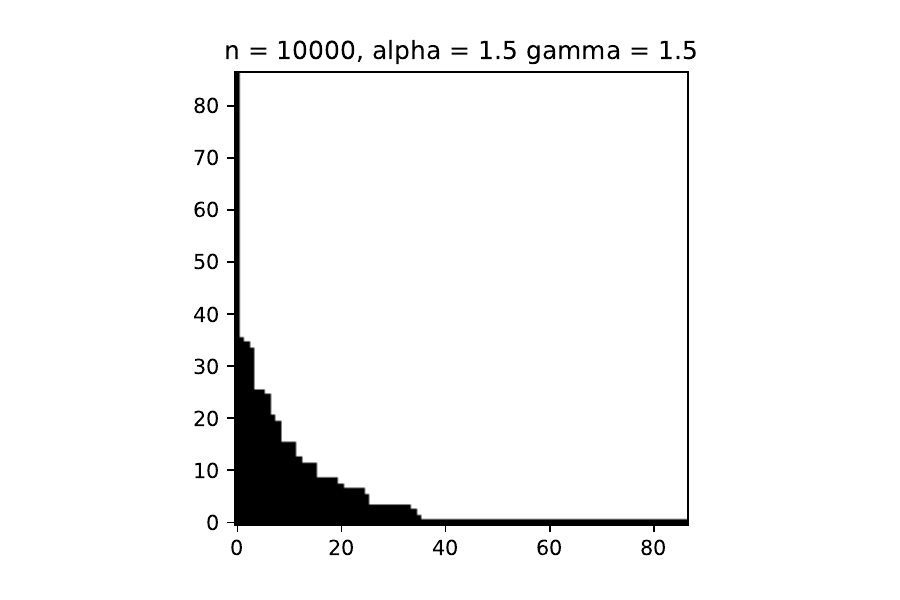}
\caption{\label{fig:1}Empirical graphons of the adapted model, with $\alpha = 1.5, \gamma = 0.5, 0.8, 1.1, 1.5$, and $n=10000$. Vertices are labeled according to decreasing vertex values.}
\end{figure}

To examine the structure of the random graph away from criticality, we introduce an auxiliary parameter $\gamma>0$ and set $a_n^\alpha=n^\gamma \log n$. This explicit parametrization has some nice implications. First, since
\begin{equation}
\ccbb{X_iX_j>a_n} = \ccbb{U_iU_j<\frac{1}{n^\gamma \log n}},
\end{equation}
$\alpha$ becomes irrelevant as a parameter, and we can concentrate on tuning the parameter $\gamma$ solely. The sub-critical regime $a_n \gg n^{2/\alpha}$ translates to $\gamma \geq 2$ and the super-critical regime $a_n \ll n^{2/\alpha}$ translates to $\gamma<2$. Second, under this parametrization, the average number of edges of the graph exhibits an asymptotic growth order that is entirely dependent on $\gamma$,
\equh\label{eq:E_n}
\esp |E_n| \sim \binom n2 \alpha a_n^{-\alpha}\log (a_n)\sim \frac{\gamma} 2 n^{2-\gamma}.
\eque
Note the difference between this order and the order of expected number of edges in the clique,
\begin{equation}
\left(\esp K_{n, 0}\right)^2 \sim \frac{n^2}{a_n^\alpha} \sim \frac{n^{2-\gamma}}{\log n}.
\end{equation}
This implies that in asymptotics, most of the edges are coming from followers. See Figure \ref{fig:1} for some simulations.

Other than the simple calculations presented above, many other distinguishing features of our adapted model could be derived explicitly. With some abuse of notation, let us denote the set of non-isolated vertices still by $V_n$. We compute the probability that a given vertex is non-isolated.
\begin{align}
\proba\pp{\max_{i=2,\dots,n}X_1X_{i}>a_n} & = \int_1^\infty \alpha x^{-\alpha-1}dx \bb{1-\proba\pp{X_1\le \frac{a_n}{x}}^{n-1}} \notag \\
& = \frac1{a_n^\alpha} + \int_1^{a_n} \alpha x^{-\alpha-1}dx \bb{1-\pp{1-\pp{\frac {x}{a_n}}^\alpha}^{n-1}} \notag \\
& = \frac 1{a_n^\alpha}+\frac n{a_n^\alpha}\int_{n/a_n^\alpha}^n y^{-2}\bb{1-\pp{1-\frac yn}^{n-1}}dy,
\end{align}
where the last step consists of a change of variables $(x/a_n)^\alpha=y/n$. In the case $\gamma\in[1,2]$, the lower limit of integration $n/a_n^\alpha \rightarrow 0$, giving
\begin{align}
\esp |V_n| & = n\proba\pp{\max_{i=2,\dots,n}X_1X_{i}>a_n} \sim \frac {n^2}{a_n^\alpha}\int_{n/a_n^\alpha}^n y^{-2}(1-e^{-y})dy \notag \\
& \sim \frac {n^{2-\gamma}}{\log n} \left(\int_{n^{1-\gamma}/\log n}^1 y^{-1}dy+\int_{1}^n y^{-2}dy\right) \notag \\
& \sim \frac {n^{2-\gamma}}{\log n} \left(\log\log n-(1-\gamma)\log n+1\right).
\end{align}
In the case $\gamma\in(0,1)$, the lower limit of integration $n/a_n^\alpha \rightarrow \infty$, giving
\begin{align}
\esp |V_n| & = n\proba\pp{\max_{i=2,\dots,n}X_1X_{i}>a_n} \sim \frac {n^2}{a_n^\alpha}\int_{n/a_n^\alpha}^n y^{-2} dy \notag \\
& \sim \frac {n^{2-\gamma}}{\log n} \int_{n^{1-\gamma}/\log n}^n y^{-2}dy \sim n.
\end{align}
The case $\gamma>2$ is clear, as $\esp |E_n| \to 0$. We summarize below the expected number of non-isolated vertices of the random graph in all parameter regimes:
\begin{equation}
\esp |V_n|\sim
\begin{cases}
0 & \gamma>2, \\
(\gamma-1)n^{2-\gamma} & \gamma\in(1,2],\\
n \frac{\log \log n}{\log n} & \gamma = 1,\\
n & \gamma\in(0,1).
\end{cases}
\end{equation}
\textit{All our calculations so far seem consistent and suggest that there is no random graph in the limit at the sub-critical regime $\gamma\geq 2$. Contrarily, at the super-critical regime $\gamma\in(0,2)$, a limit random graph exists.} We will respectively investigate these asymptotic phenomena in detail in Sections \ref{sec:sub_critical} and \ref{sec:super_critical}.

Letting $A_{n,i}$ denote the event that vertex $i$ is not isolated and $B_{n,(i,j)}$ the event that $(i,j)$ is an edge, the next quantity we will compute is $\rho_{n,(i,j)} = \proba(B_{n,(i,j)}\mid A_{n,i}\cap A_{n,j})$. If this conditional probability has a limit, say $\rho_{(i,j)}$, then it may be interpreted as the probability of having an edge between two non-isolated vertices of the limit graph.
\begin{multline}
\PR\left(A_{n, 1} \cap A_{n, 2} \cap B_{n, (1, 2)}^c \right)
=2\int_1^{\infty} \alpha x_1^{-\alpha-1} dx_1 \int_1^{x_1} \alpha x_2^{-\alpha-1} dx_2 \ \inddd{x_1x_2<a_n} \left(1-\PR\left(X_1\leq \frac{a_n}{x_2}\right)^{n-2}\right) \\
=2\int_1^{\sqrt{a_n}} \alpha x_2^{-\alpha-1} dx_2 \int_{x_2}^{\frac{a_n}{x_2}} \alpha x_1^{-\alpha-1} dx_1 \left(1-\left(1-\left(\frac{x_2}{a_n}\right)^{\alpha}\right)^{n-2}\right) \\
\sim \frac{2n^2}{a_n^{2\alpha}} \int_{\frac{n}{a_n^\alpha}}^{\frac{n}{a_n^{\alpha/2}}} y_2^{-2} \left(1-e^{-y_2}\right) dy_2 \int_{y_2}^{\frac{n^2}{a_n^\alpha}\frac{1}{y_2}} y_1^{-2} dy_1,
\end{multline}
where for the first equality, w.l.o.g. we assumed that $x_1>x_2$. The indicator function constraint then gives $x_2<\sqrt{a_n}$ in the second equality. We then apply a change of variables $(x_i/a_n)^{\alpha}=y_i/n$ for $i=1,2$. As was explained in the above calculation for $\esp |V_n|$, a standard asymptotic study yields
\begin{equation}
\PR\left(A_{n, 1} \cap A_{n, 2} \cap B_{n, (1, 2)}^c \right)\sim
\begin{cases}
\frac{2}{n^{\gamma-1}\log n} & \gamma \geq 1,\\
1 & \gamma\in(0,1).
\end{cases}
\end{equation}
Complementarily,
\begin{multline}
\PR\left(A_{n, 1} \cap A_{n, 2} \cap B_{n, (1, 2)}\right)
=\int_1^\infty \alpha x_1^{-\alpha-1} dx_1 \int_1^\infty \alpha x_2^{-\alpha-1} dx_2 \ \inddd{x_1x_2>a_n} \\
=\int_1^{a_n} \alpha x_1^{-\alpha-1} dx_1 \int_{\frac{a_n}{x_1}}^\infty \alpha x_2^{-\alpha-1} dx_2 + \int_{a_n}^\infty \alpha x_1^{-\alpha-1} dx_1 \int_1^\infty \alpha x_2^{-\alpha-1} dx_2
\sim \frac{\gamma}{n^\gamma} \hspace{.2cm} \text{ for all } \gamma>0.
\end{multline}
Combining the above results, we have $\rho_{(i,j)}=0$ for all $\gamma>0$. The zero conditional probability of an edge connecting two non-isolated vertices may look puzzling at first sight when $\gamma \in (0,2)$ since we expect the existence of an infinite random graph in this parameter region. A possible interpretation is that, as the number of non-isolated vertices grows infinite in an appropriate sense, given any two vertices the probability of seeing an edge between them is zero.

\section{Sub-critical regime}\label{sec:sub_critical}
Since there is no random graph in distribution in the limit, the sub-critical regime $a_n\gg n^{2/\alpha}$ is relatively not as interesting as compared to the super-critical regime $a_n\ll n^{2/\alpha}$. Nevertheless, the limit object for the sub-critical regime captures some intriguing characteristics, as we will describe in this section. For explicitness, we take the i.i.d.~$X_i$ to have pdf $\alpha x^{-\alpha-1}dx\inddd{x\ge 1}$ and set $a_n^\alpha=n^\gamma \log n$ with $\gamma \geq 2$.

Recall that $K_{n,0} = \summ i1n \inddd{X_i/\sqrt{a_n}>1}$ denotes the number of vertices with large weight (vertex value $>\sqrt{a_n}$) and is Binomial distributed with parameter $\left(n, a_n^{-\alpha/2}\right)$. These vertices are referred to as clique vertices if they are in addition non-isolated. Since
\begin{equation}
K_{n, 0} \rightarrow \esp K_{n, 0} \sim \frac{n^{1-\gamma/2}}{\log ^{1/2}n} \to 0
\end{equation}
in probability in the sub-critical regime, only a conditional limit theorem is worth investigating in this case. The conditional law of $\calL(K_{n, 0}\left|\right.K_{n, 0}\geq 1)$ is easy to derive, with
\begin{equation}\label{fake1}
\PR\left(K_{n, 0}\geq 1\right)=1-\left(1-\frac{1}{a_n^{\alpha/2}}\right)^n \sim \frac{n}{a_n^{\alpha/2}}=\frac{n^{1-{\gamma/2}}}{\log ^{1/2}n},
\end{equation}
\begin{equation}\label{fake2}
\PR\left(K_{n, 0}= 1\right)=n\left(\frac{1}{a_n^{\alpha/2}}\right) \left(1-\frac{1}{a_n^{\alpha/2}}\right)^{n-1} \sim \frac{n}{a_n^{\alpha/2}}=\frac{n^{1-{\gamma/2}}}{\log ^{1/2}n}.
\end{equation}
We conclude that given the appearance of a non-trivial random graph, the clique part (conditioning on non-empty) typically only contains one vertex. There is a fine point when computing the probability that $K_{n,0}=1$ though. This event does not necessarily imply the appearance of a star graph as the edge number may still be zero; just having one large weight vertex is not enough. We demonstrate this subtlety below.
\begin{multline}
\PR\left(\text{one clique vertex}\right)\\
=n\int_{\sqrt{a_n}}^{a_n} \alpha x_1^{-\alpha-1} dx_1 \left[\left(1-\frac{1}{a_n^{\alpha/2}}\right)^{n-1}-\left(1-\left(\frac{x_1}{a_n}\right)^\alpha\right)^{n-1}\right]\\
+n\int_{a_n}^{\infty} \alpha x_1^{-\alpha-1} dx_1 \left(1-\frac{1}{a_n^{\alpha/2}}\right)^{n-1}.
\end{multline}
Here we eliminate the situation where $K_{n, 0}=1$, but no edge is formed between the clique vertex $X_1$ and the follower vertices $X_2, \dots, X_n$. The scalar $n$ indicates that the clique could be centered at any vertex. The second term on the right is of order $na_n^{-\alpha} \ll n^{1-\gamma}$, while the first term, after a change of variables $(x_1/a_n)^{\alpha}=y_1/n$, asymptotically becomes
\begin{equation}
\frac{n^2}{a_n^\alpha} \int_{\frac{n}{a_n^{\alpha/2}}}^n y_1^{-2} \left(1-e^{-y_1}\right) dy_1 \sim \frac{n^2}{a_n^\alpha} \log \left(\frac{a_n^{\alpha/2}}{n}\right) \sim \left(\frac{\gamma}{2}-1\right)n^{2-\gamma}.
\end{equation}
Combining the above analysis, the correct probability of a star graph with a lone vertex in the clique behaves like $(\gamma/2-1)n^{2-\gamma}$, which is smaller than (\ref{fake1}) (\ref{fake2}).

Let us delve deeper into the structure of this star graph. Notice that $K_{n, 1}$ by itself counts the number of followers in this case; $K_{n, j}=0$ for all $j>1$ automatically by construction. We have
\begin{align}
\PR\left(\text{one clique vertex}, K_{n, 1}=1\right)
&=n(n-1) \int_{\sqrt{a_n}}^{a_n} \alpha x_1^{-\alpha-1} dx_1 \int_{\frac{a_n}{x_1}}^{\sqrt{a_n}} \alpha x_2^{-\alpha-1} dx_2 \left(1-\left(\frac{x_1}{a_n}\right)^\alpha\right)^{n-2}\notag \\
&\sim \frac{n^2}{a_n^\alpha} \int_{\frac{n}{a_n^{\alpha/2}}}^n y_1^{-1} e^{-y_1} dy_1 \sim \frac{n^2}{a_n^\alpha} \log \left(\frac{a_n^{\alpha/2}}{n}\right) \sim \left(\frac{\gamma}{2}-1\right)n^{2-\gamma}.
\end{align}
Here the scalers $n$ and $n-1$ in the equality indicate that the clique could be centered at any vertex and the follower could come from any of the remaining vertices. A standard asymptotic analysis then yields the asymptotic order after a change of variables $(x_1/a_n)^{\alpha}=y_1/n$. This probability is asymptotically the same as having a lone clique star graph that was established previously. We state this finding.

\begin{theorem}
Consider the adapted model at the sub-critical regime ($a_n^\alpha=n^\gamma \log n$ with $\gamma \geq 2$). Given that the graph is non-empty, in the limit predominantly it has exactly two vertices, one clique vertex and one follower vertex.
\end{theorem}

A physical interpretation of this phenomenon might be the following: For a typical behavior, with probability going to one we would not see any graph eventually. In the rare event that we do see one, we would need certain `extra energy' (than typical) to push some of the $X_i$ values up, and the most `economical' way to do so is to push one up to the clique and another up as a follower. Pushing up two to the clique or pushing up more than one follower or any other construction, by comparison, might be too costly.

\section{Super-critical regime}\label{sec:super_critical}
In this section we will examine the structure of the adapted model in the more intriguing super-critical regime $a_n \ll n^{2/\alpha}$, where the limit random graph is expected to take a non-trivial form. Recall that $X_i\eqd U_i^{-1/\alpha}$ and $U_i$ are i.i.d.~$(0,1)$-uniform random variables. Denote the tail distribution by
\equh\label{eq:RV}
\wb F(x) = \proba(X_1>x) = (x\vee 1)^{-\alpha}.
\eque In a complementary manner, we also write $F(x) = \proba(X_1\le x)$ for the cumulative distribution function. For ease of notation, let $K_n:=K_{n,0} = \summ i1n \inddd{X_i/\sqrt{a_n}>1}$ denote the number of vertices in the clique. In the super-critical regime, $a_n= o(n^{2/\alpha})$ gives
\equh\label{eq:super_critical}
\sigma_n^2:=\esp K_n = n\wb F(\sqrt{a_n}) \to \infty.
\eque
Since $K_n$ is Binomial distributed with parameter $\left(n, \wb F(\sqrt{a_n})\right)$, $K_n \sim \sigma_n^2$ in probability, and the random variable $K_n$ is well-concentrated around $\sigma_n^2$.

Introduce two i.i.d.~sequences of random variables $\{Y_{n,i}\}_{i\in\N}$ and $\{Z_{n,i}\}_{i\in\N}$ with
 \begin{align}\label{eq:Yn1}
 \proba(Y_{n,1}>y)& = \frac{\wb F(y\sqrt{a_n})}{\wb F(\sqrt {a_n})}=y^{-\alpha},\mmas n\to\infty \mfa y>1,\\
 \proba\pp{Z_{n,1}\le x} & = \frac{F(x)}{F(\sqrt{a_n})}, \mbox{ } x\in(0,\sqrt{a_n}].
 \end{align}
In other words, $\{Y_{n,i}\}_{i\in\N}$ are i.i.d.~with law as $\calL(X_1\left|\right.X_1>\sqrt{a_n})$ (with scaling adjustment) and $\{Z_{n,i}\}_{i\in\N}$ are i.i.d.~with law as $\calL(X_1\left|\right.X_1\leq \sqrt{a_n})$. Assume further that these two sequences are independent.
Then for every $n\in\N$, given $K_n$, the values of $\{X_i\}_{i=1,\dots,n}$ corresponding to those larger than (less than resp.) the threshold $\sqrt{a_n}$ share the same joint law of $Y_{n,1},\dots,Y_{n,K_n}$ ($Z_{n,1},\dots,Z_{n,n-K_n}$ resp.).
We order $\{Y_{n,i}\}_{i=1,\dots,K_n}$ in {\em increasing} order statistics
 \begin{equation}
 Y_{n,K_n:K_n}>\cdots>Y_{n,1:K_n}>\sqrt {a_n}>\frac{a_n}{Y_{n,1:K_n}}>\cdots>\frac{a_n}{Y_{n,K_n:K_n}},
 \end{equation}
 where listed on the right hand side are the thresholds for different groups of followers.

Define the statistics
\begin{equation}\label{tau_n}
 \tau_n(x):=\frac{a_n}{Y_{n,\ceil{x K_n}:K_n}}, \hspace{.2cm} x\in(0,1).
\end{equation}
 We are interested in the asymptotic behavior of the \emph{height function}
 \begin{equation}\label{eq:H_n}
 H_n(x) := \summ i1{n-K_n}\inddd{Z_{n,i}>\tau_n(x)}.
 \end{equation}
This construction associates the law of $H_n(x)$ to that of the number of not-in-clique vertices that are connected to the vertices in the clique corresponding to those top $\ceil{xK_n}$-values of $\{Y_{n,i}\}_{i=1,\dots,K_n}$. At one end, $H_n(1)$ is the number of followers of the leader from the clique, i.e., $n-K_n$, thus $H_n(1) \sim n$ with high probability. At the other end, we take $H_n(0)\equiv 0$ by convention.
For notational convenience, set
\begin{equation}
B_{n,i}(x) := \inddd{Z_{n,i}>\tau_n(x)}.
\end{equation}

We are now ready to state our main result, which says that the limit fluctuation of the height function has two independent components, one as a generalized Brownian bridge, the other as a time-changed Brownian motion.

\begin{theorem}\label{thm:1}
Consider the adapted model at the super-critical regime ($a_n^\alpha=n^\gamma \log n$ with $\gamma<2$). Let $\sigma_n^2$ and $H_n(x)$ be defined as in \eqref{eq:super_critical} and \eqref{eq:H_n}. We have
\equh\label{eq:joint_fdd}
\frac1{\sigma_n}\ccbb{H_n(x) -{\sigma_n^2}\frac{x}{1-x} }_{x\in[0,1)} \stackrel{f.d.d.}\weakto\ccbb{\B_{x/(1-x)}+\G_x}_{x\in[0,1)},
\eque
where $f.d.d.$ indicates convergence of finite-dimensional distributions,
$\{\B_t\}_{t\in[0,\infty)}$ is a standard Brownian motion,
$\{\G_x\}_{x\in[0,1)}$ is a generalized Brownian bridge with covariance function
\begin{equation}
\cov\pp{\G_x,\G_y} =
\frac{\min(x,y)(1-\max(x,y))}{(1-x)^2(1-y)^2}, \hspace{.2cm} x,y\in[0,1),
\end{equation}and $\B$ and $\G$ are independent.
\end{theorem}

Note that throughout, the index $x$ is strictly less than $1$ (covariance explodes as $x\uparrow 1$). On the other hand, convergence at $x=0$ is clear, as both sides equal zero.

\vskip.1truein

\noindent \textbf{Implications of Theorem \ref{thm:1}.} Recall that $\sigma_n^2=\esp K_n$ denotes the average number of vertices in the clique. From Theorem \ref{thm:1}, we may deduce that
\begin{equation*}
\esp H_n(x) \sim \sigma_n^2 \frac{x}{1-x},
\end{equation*}
\begin{equation}
\var \pp{H_n(x)} \sim \sigma_n^2 \frac{x}{1-x} \left(1+\frac{1}{(1-x)^2}\right).
\end{equation}
Since we established the theorem using increasing order statistics, this introduces a simple transformation $x \mapsto 1-x$ to the simulations in Figure \ref{fig:1}. Then for $x\in (0, 1]$,
\begin{equation}\label{eq:h}
h(x):=1+\limn \frac{\esp H_n(x)}{\sigma_n^2}=\frac{1}{x}
\end{equation}
should give the boundary line in question, where the extra $1$ in the above expression comes from the clique-clique contribution. This is a universal result independent of the parameters. Having the same asymptotic order for the expected value and the variance of the height function $H_n(x)$ also explains why the simulations look so regular.

Let $W_n$ denote the graphon of our model with $n$ vertices without scaling (a $\{0,1\}$-valued function on $[0,n]^2$). Following explanations above,
\begin{equation}
W_n'(x,y) = W_n\left(\esp K_n \cdot x, \hspace{.1cm} \esp K_n \cdot y\right)
\end{equation}
has the deterministic limit $W(x,y)=\inddd{xy\le 1}, x,y\in(0,\infty)$, with boundary line $y=1/x$ (\ref{eq:h}). The asymptotics is confirmed by the simulations in Figure \ref{fig:1}. This mode of convergence feels very close to the stretched convergence mode \eqref{eq:stretched}. Our stretching acts in a similar way as in \eqref{eq:stretched}, but by a different stretching order, as $\nn {W_n}_1 = \esp |E_n| \sim (\gamma/2)n^{2-\gamma}$ while $\esp K_n=na_n^{-\alpha/2}\sim n^{1-\gamma/2}/(\log n)^{1/2}$. So there is an extra $\log$ term in our stretching as compared to \eqref{eq:stretched}. Furthermore, contrary to the stretched convergence mode, the limit graphon $W$ is not $L^1$-integrable. Our adapted model may therefore be viewed as an example that lies at the boundary between rescaled convergence and stretched convergence.

\vskip.1truein

To prove Theorem \ref{thm:1}, we shall decompose $H_n$ further and identify the relevant contribution from the different parts of the statistics to $\B$ and $\G$. For this purpose, we introduce $\calK_n := \sigma(K_n,Y_{n,1},\dots,Y_{n,K_n})$, and
\equh
\what p_n(x) := \proba\pp{Z_{n,i}>\tau_n(x)\mmid \calK_n}
 \qmand
\label{eq:pn}
p_n(x):= \frac{\wb F(\sqrt{a_n})}{F(\sqrt {a_n})}\pp{\theta_n(1-x)-1},
\eque
where
\begin{equation}\label{eq:theta_n}
\theta_n(x) = \begin{cases}
x\inv & x>\wb F(\sqrt {a_n}) = a_n^{-\alpha/2},\\
a_n^{\alpha/2} & x\le a_n^{-\alpha/2}.
\end{cases}
\end{equation}
Write
\begin{align}
\wb H_n(x) & =  \summ i1{n-K_n}\pp{B_{n,i}(x)-p_n(x)} \notag \\
& =\summ i1{n-K_n}(B_{n,i}(x)-\what p_n(x)) + (n-K_n)(\what p_n(x)-p_n(x))=: \wb H_{1,n}(x)+\wb H_{2,n}(x).
\end{align}
We will show that $\sigma_n\inv\wb H_{1,n}$ and $\sigma_n\inv \wb H_{2,n}$ converge to $\B$ and $\G$, respectively. %The convergence is actually stronger than finite-dimensional distributions.
We first examine $\wb H_{2,n}$.
\begin{lemma}\label{lem:pn_hat}
With the notations above, under the assumptions in Theorem \ref{thm:1}, we have
\begin{equation}
\frac n{\sigma_n}\ccbb{\what p_n(x)-p_n(x)}_{x\in(0,1)}\stackrel{f.d.d.}\weakto \ccbb{\G_x}_{x\in(0,1)}.
\end{equation}
\end{lemma}
\begin{proof}
We begin the analysis of asymptotics by first examining the i.i.d.~$Y_{n,i}$. Definition \eqref{eq:Yn1} implies that $Y_{n,1}\eqd U^{-1/\alpha}$, where $U$ is a uniform random variable on $(0,1)$. Set
$W_{n,i}:=Y_{n,i}\inv \eqd U^{1/\alpha}$
and $W_n:=W_{n,1}$. % and $V_n = V_{n,1}$. %(\alert{TODO} to write all the proof in $W_n$ only.)
We need some background on quantile processes from Shorack \cite{shorack72functions} \cite{shorack73convergence}. Following notations from earlier, let $F_Z$ be the cumulative distribution function of a random variable $Z$ and $\wb F_Z$ its tail probability function. Let $F_Z\inv$ denote the left-continuous inverse function of $F_Z$. Let $\mathbb F\inv _{Z,n}(x)$ denote the quantile process of i.i.d.~copies $Z_{1},\dots,Z_{n}$. It follows that $\F\inv_{Z,n}(x) = Z_{\floor{xn}:n}$, where $Z_{i:n}$ is the $i$th smallest order statistic of $Z_1,\dots,Z_n$.
 Then $F_{W_n}(W_{n})$ is a uniform random variable, and for all $m_n\to\infty$,
\equh\label{eq:quantile_convergence}
\sqrt {m_n}\ccbb{\mathbb Q_{n,m_n}(x)-x}_{x\in[0,1]}:=\sqrt {m_n}\ccbb{F_{W_n}\circ\mathbb F_{W_n,m_n}\inv(x)-x}_{x\in[0,1]}\weakto \ccbb{\B^{br}_x}_{x\in[0,1]}.
\eque
Here, $\{\B^{br}_x\}_{x\in[0,1]}$ is a standard Brownian bridge, a centered Gaussian process with
\begin{equation}
\cov\pp{\B^{br}_x,\B^{br}_y} = \min(x,y)(1-\max(x,y)), \hspace{.2cm} x,y\in[0,1]
\end{equation}
in $D([0,1])$. $\mathbb Q_{n,m_n}(x)$ so defined has the law of the quantile process of $m_n$ i.i.d.~uniform random variables. Furthermore, from (\ref{tau_n}),
\begin{equation}
\tau_n(x)= \sqrt{a_n}\F\inv_{W_n,K_n}(1-x).
\end{equation}

Combining the above observations, we have
\begin{equation}
\what p_n(x) = \frac{\wb F(a_n/Y_{n,\ceil{xK_n}:K_n})-\wb F(\sqrt{a_n})}{F(\sqrt{a_n})} = \what \rho_n(x) \frac{\wb F(\sqrt{a_n})}{F(\sqrt{a_n})},
\end{equation}
with
\begin{equation}
\what \rho_n(x)  =  \frac{\wb F(\sqrt{a_n}\mathbb F_{W_n,K_n}\inv(1-x))}{\wb F(\sqrt {a_n})}-1.
\end{equation}
%Note that $F_{V_n}\inv(x) = 1/F_{W_n}\inv(1-x)$, and their order statistics have a corresponding relation.
We rewrite $\what \rho_n$ as a function of $\mathbb Q_n:=\mathbb Q_{n,K_n}$. Namely,
\begin{align}
\what \rho_n(x) & %=  \frac{\wb F(a_n/\mathbb F_{Y_n,K_n}\inv(x))}{\wb F(\sqrt {a_n})}-1
%= \frac{\wb F(\sqrt{a_n}\mathbb F_{W_n,K_n}\inv(1-x))}{\wb F(\sqrt {a_n})}-1
 = \frac{\wb F(\sqrt{a_n}F_{W_n}\inv\circ F_{W_n}\circ \mathbb F_{W_n, K_n}\inv(1-x))}{\wb F(\sqrt {a_n})}-1 \notag \\
 &
 = \frac{\wb F(\sqrt{a_n}F_{W_n}\inv\circ \mathbb Q_{n}(1-x))}{\wb F(\sqrt {a_n})}-1.
\end{align}

Let us make these calculations more explicit. For $x\in(0,1)$, $F_{W_n}(x) = x^\alpha$ and $F_{W_n}\inv(x) = x^{1/\alpha}$. Hence
\begin{equation}
\frac{\wb F\pp{\sqrt{a_n}F_{W_n}\inv(x)}}{\wb F(\sqrt{a_n})}=\theta_n(x),
\end{equation}
where $\theta_n(x)$ is defined as in \eqref{eq:theta_n}. This gives
\begin{equation}
\what \rho_n(x) = \theta_n(\mathbb Q_{n}(1-x))-1.
\end{equation}
We have
\equh\label{eq:delta}
\frac n{\sigma_n}(\what p_n(x) - p_n(x)) = \frac1{F(\sqrt{a_n})}\cdot \frac{\sigma_n}{\sqrt{K_n}}\cdot \sqrt{K_n}\pp{\theta_n(\mathbb Q_n(1-x))-\theta_n(1-x)}.
\eque
Note that $\theta_n'(x) = -x^{-2}$ for all $x>0$ for $n$ large enough. A standard application of the delta method applied to \eqref{eq:quantile_convergence} and \eqref{eq:delta} then yields the desired weak convergence.
\end{proof}

The examination of $\wb H_{1,n}$ comes next. Following \cite{shorack73convergence} we can actually construct, on a different probability space, for each $n\in\N$ copies $\{\wt Y_{n,i}, \wt K_n\}_{i\in\N}$ of $\{Y_{n,i},K_n\}_{i\in\N}$,  such that $\wt K_n \sim  \sigma_n^2$ almost surely, and the convergence \eqref{eq:quantile_convergence} is in the almost sure sense.
 Note that this coupling construction does not necessarily imply  that the joint laws of $(\wt Y_{n,i},\wt Y_{m,j})$ are the same as $(Y_{n,i},Y_{m,j})$ for $m\ne n$, but we will not need such joint laws in the sequel. For ease of notation, we still let $Y_{n,i},W_{n,i},K_n,\tau_n(x)$ denote $\wt Y_{n,i},\wt W_{n,i},\wt K_n,\wt\tau_n(x)$ respectively, and emphasize that we are working on this different probability space simply by saying {\em for the coupled model}. Recall that $\calK_n := \sigma(K_n,Y_{n,1},\dots,Y_{n,K_n})$. Further define $\calK:=\sigma(\{\calK_n\}_{n\in\N})$ for the coupled model. We also continue to assume $Z_{n,i}$ as before, independent from all other random variables discussed so far in this paragraph.

We start by noticing that for each $n\in\N$, given $\calK_n$, $\{B_{n,i}(x)\}_{i\in\N}$ are i.i.d.~Bernoulli random variables with parameter $\what p_n(x)$ for every $x$. Moreover, they are nested in the sense that $B_{n,i}(x) = 1$ implies $B_{n,i}(y) = 1$ for all $y\in(x,1)$. We fix $x\in(0,1)$ (there is nothing to prove for $x = 0$). Then, conditioning on $\calK_n$, for each $n\in\N$,
\begin{equation}
\wb H_{1,n}(x) = \summ i1{n-K_n}(B_{n,i}(x)-\what p_n(x))
\end{equation}
is a partial sum of (centralized) i.i.d.~Bernoulli random variables with parameter $\what p_n(x)$. Now, from \eqref{eq:quantile_convergence}, we know that
\begin{equation}
\mathbb Q_{n}(x):=F_{W_n}\circ\F_{W_n,K_n}\inv(x)\to x \mbox{ almost surely,}
\end{equation}
and $\limn F_{W_n}(x) =x^\alpha$ for $x\in (0,1)$.
Therefore,
the conditional variance is
\begin{align}
(n-K_n)\what p_n(x)(1-\what p_n(x)) & \sim n\wb F(\sqrt{a_n})\what\rho_n(x) = n\wb F(\sqrt{a_n})(\theta_n(\Q_n(1-x))-1) \notag \\%= \pp{\frac{\wb F(\sqrt {a_n}\F_{W_n,K_n}\inv(1-x))}{\wb F(\sqrt{a_n})}-1}\\
& %\sim\sigma_n^2\pp{\frac{\wb F(\sqrt{a_n} (1-x)^{1/\alpha})}{\wb F(\sqrt{a_n})}-1}
 \sim \sigma_n^2\pp{\frac1{1-x}-1}=\sigma_n^2 \frac{x}{1-x}
\end{align}
almost surely, where in the $\sim$ step we used the coupling $\Q_n(x)\to x$ almost surely and $\limn\theta_n(x)= x\inv$.
Then, by the central limit theorem for triangular arrays of i.i.d.~random variables, we have that
\begin{equation}
\calL\pp{\frac1{\sigma_n}\wb H_{1,n}(x)\mmid \calK}\to \calL\pp{\calN\pp{0,\frac x{1-x}}} \mbox{ a.s.}
\end{equation}
The above statement is interpreted as {\em almost-sure weak convergence}, meaning that
\begin{multline}\label{eq:asweakconvergence}
\limn \esp\bb{\phi\pp{\frac1{\sigma_n}\wb H_{1,n}(x)}\mmid\calK} = \esp \phi\pp{\pp{\frac x{1-x}}^{1/2}Z} \\
\mbox{ almost surely for all continuous and bounded functions $\phi:\R\to\R$,}
\end{multline}
where $Z$ on the right hand side is a standard Gaussian random variable. We shall use $\calL_{f.d.d.}(\{\sigma_n\inv\wb H_{1,n}(x)\}_{x\in[0,1)}\mid\calK)$ in Lemma \ref{lem:aswc} below for the corresponding {\em almost-sure weak convergence of finite-dimensional distributions} of $\wb H_{1,n}$.

This argument can be readily extended to the multivariate central limit theorem, and it suffices to compute the covariance. Alternatively, by a standard Poissonization argument one sees immediately that the limit Gaussian process has independent increments. So, the limits of finite-dimensional distributions of $\{\sigma_n\inv\wb H_{1,n}(x)\}_{x\in[0,1)}$ are the corresponding ones of $\B_{x/(1-x)}$.
We have thus proved the  following.
\begin{lemma}\label{lem:aswc}For the coupled model, we have
\begin{equation}
\calL_{f.d.d.}\pp{\frac1{\sigma_n}\{\wb H_{1,n}(x)\}_{x\in[0,1)}\mmid \calK} \to \calL_{f.d.d.}\pp{\{\B_{x/(1-x)}\}_{x\in[0,1)}} \mbox{ a.s.},
\end{equation}
where the interpretation of the above expression is explained right after \eqref{eq:asweakconvergence}.
\end{lemma}

\begin{proof}[Proof of Theorem \ref{thm:1}]
Combining Lemmas \ref{lem:pn_hat} and \ref{lem:aswc} we obtain immediately that
\begin{equation}
\frac1{\sigma_n}\ccbb{\wb H_n(x)}_{x\in[0,1)}\stackrel{f.d.d.}\weakto \ccbb{\B_{x/(1-x)}+\G_x}_{x\in[0,1)}.
\end{equation}
However, $\sigma_n\inv\wb H_n(x)$ has a slightly different centering from the one in Theorem \ref{thm:1}. Using \eqref{eq:super_critical}, \eqref{eq:pn} and \eqref{eq:theta_n}, the difference is
\begin{align}
&\frac1{\sigma_n}\pp{\sigma_n^2 \frac{x}{1-x} - (n-K_n)p_n(x)} \notag \\
=&\frac1{\sigma_n}\pp{\sigma_n^2 \frac{x}{1-x} - \sigma_n^2(\theta_n(1-x)-1)+\sigma_n^2(\theta_n(1-x)-1)- (n-K_n)p_n(x)} \notag \\
= & \sigma_n \left((1-x)^{-1}-\theta_n(1-x)\right) + \frac1{\sigma_n}\pp{(n-\esp K_n)p_n(x) - (n-K_n)p_n(x)} \notag \\
= & \sigma_n \left((1-x)^{-1}-\theta_n(1-x)\right) + \frac{K_n-\esp K_n}{\sigma_n}p_n(x),
\end{align}
which converges to zero in $L^2$.\end{proof}

\section*{Acknowledgements}
Mei Yin's research was supported in part by the University of Denver's Faculty Research Fund 84688-145601. She acknowledges helpful conversations with Yufei Zhao, and is particularly grateful to Yizao Wang for many constructive comments. She also thanks the anonymous reviewer for their insightful comments and suggestions.

\end{document}